\documentclass{amsart}
\usepackage[T1]{fontenc}
\usepackage[applemac]{inputenc} 
\usepackage[american, frenchb]{babel}

\def\ladate{{\small version du 10 mars 2003}}

\newtheorem{theoreme}{Théorème}[]

\theoremstyle{definition}

\newcommand\NN{{\mathbb N}}

\newcommand\RR{{\mathbb R}}

\newcommand\cA{{\mathcal A}}
\newcommand\cB{{\mathcal B}}

\newcommand\cF{{\mathcal F}}

\renewcommand{\Re}{{\rm Re}}

\begin{document}

%%\textbf{Analyse Harmonique/\textit{Harmonic Analysis}}
\rightline{\ladate}
\rightline{math/0302102}

\bigskip

\title[Des \'equations de Dirac-Schr\"odinger pour
Fourier]{Des \'equations de Dirac et de Schr\"odinger pour
la transformation de Fourier}

\author{Jean-Fran\c cois Burnol}

\thanks{%%\textbf{Note présentée par J.-P. Kahane}\\
La version datée 10 mars 2003 apporte à celle du 11 février
les équations additionnelles  \eqref{eq:phidet},
\eqref{eq:psidet},  \eqref{eq:schrodbis++}, \eqref{eq:schrodbis+-},
\eqref{eq:schrodbis-+}, \eqref{eq:schrodbis--}, ainsi que
\eqref{eq:iso} et \eqref{eq:jost} et s'accompagne d'une
discussion plus précise de ceux des résultats de l'auteur qui feront
l'objet d'une rédaction plus détaillée par ailleurs.}

\begin{abstract}
Dyson a associé aux déterminants de Fredholm des noyaux de Dirichlet
pairs (resp. impairs) une équation de
Schrödinger sur un demi-axe et a employé les méthodes du
scattering inverse de Gel'fand-Levitan et de Marchenko, en
tandem, pour étudier l'asymptotique de ces
déterminants. Nous avons proposé suite à notre mise-au-jour
de l'opérateur conducteur de chercher à réaliser la
transformation de Fourier elle-même comme un scattering, et
nous obtenons ici dans ce but deux systèmes de Dirac sur
l'axe réel tout entier  et qui sont associés
intrinsèquement, respectivement, aux transformations  en
cosinus et en sinus.
\end{abstract}

\maketitle

\selectlanguage{american} \title{\textbf{Dirac and
Schrödinger equations for the Fourier Transform}}

\begin{abstract}
Dyson has associated with the Fredholm determinants of the
even (resp. odd) Dirichlet kernels a Schrödinger
equation on the half-axis and has used, in tandem, the
Gel'fand-Levitan and Marchenko methods of inverse scattering
theory to study the asymptotics of these determinants. We
have proposed following our unearthing of the conductor
operator to seek to realize the Fourier transform itself as
a scattering, and we obtain here to this end two Dirac
systems on the entire real axis which are intrinsically
associated, respectively, to the cosine and to the sine
transforms.
\end{abstract}

\maketitle

\selectlanguage{french}
\bigskip

\begin{small}
Jean-Fran\c cois Burnol, Université Lille 1, UFR de Mathématiques, 
Cité Scientifique M2, F-59655 Villeneuve d'Ascq Cedex, France\\
\texttt{burnol@agat.univ-lille1.fr}
\end{small}

\centerline{\hbox to4cm{\leaders \hrule height1pt \hfill}}

\selectlanguage{american}
\section{Abridged English Version}

Dyson \cite{dyson} has associated with the Dirichlet kernels
\begin{equation}
D_s^\pm(x,y) = \frac{\sin(s(x-y))}{\pi(x-y)} \pm
\frac{\sin(s(x+y))}{\pi(x+y)}
\end{equation}
on $(0,1)$ two Schrödinger equations on 
$(0,+\infty)$ whose potentials are given as:
\begin{equation}
W_\pm(s) = -2 \frac{d^2}{ds^2}\log D_\pm(s) - 1
\end{equation}
In this definition the $D_\pm(s)$ are the Fredholm Determinants
$\det(\delta(x-y)-D_s^\pm(x,y))$ associated with the
(trace-class) kernels above. After establishing a general
link between the methods of Gel'fand-Levitan on one hand,
of Marchenko on the other hand, and the study of
Fredholm determinants, Dyson obtained in this particular case
the existence of an asymptotic development as  $s\to\infty$
of
$\log D_\pm(s)$, thus extending the earlier result of
des~Cloizeaux
and Mehta (\cite{dcmehta}) which is:
\begin{equation}
\log D_\pm(s) = -\frac{s^2}4\mp\frac12 s- \frac18\log s +
C_\pm + o(1),
\end{equation} and another heuristic argument enabled Dyson
to give the values for the constants $C_\pm$. Except for
that last point, still in need of rigorous proof,
des~Cloizeaux-Mehta's  and Dyson's results have been later
rigorously established by Deift-Its-Zhou \cite{diz} in the
framework of a Riemann-Hilbert approach (see also Widom
\cite{widom} for a proof of $\log
D_+(s)D_-(s) = -\frac{s^2}2(1+o(1))$).

Our discovery  \cite{jfprop} of the local operator theoretic
content of the explicit formulae of Number Theory led us to
propose to realize the Fourier transform itself as a
scattering \cite{jfscat}. The multiplicative Fourier
analysis of the additive Fourier transform  leads to two
unit modulus functions on the critical line $\chi_+(s)$ and
$\chi_-(s)$, the first corresponding to the cosine
transform, and the second to the sine transform. But their
phases are of the order $\gamma\log|\gamma|$ ($s =
\frac12+i\gamma$) and this escapes the framework of
scattering ``matrices'' for potentials which satisfy  the usual
hypotheses of the Gel'fand-Levitan and Marchenko methods.

We associate in this Note to the cosine transform (resp.
sine transform) a differential system which, as in
\cite{dyson}, involves the Fredholm determinants
$D_\pm(s)$. To express it as a Dirac system, we shall use,
not the variable $s=2\pi a^2\in (0,\infty)$, but
rather the variable $u = \log(a) \in (-\infty,+\infty)$. Let
us thus define $D_a^\pm$, $a>0$ to be the Dirichlet kernel
as above for $s=2\pi a^2$ and let  $d_a^\pm = \det(1 -
D_a^\pm)$ and $d_\pm(u) = d_a^\pm, a= \exp(u), u\in\RR$. We
propose the following Dirac-type system:

\begin{subequations}\label{eq:dirac}
\begin{align}
\frac{d}{du} \alpha(u) &= - \sqrt{-\frac{d^2}{du^2}\log
d_\pm(u)}\cdot\alpha(u) - \gamma \beta(u)\\
\frac{d}{du} \beta(u) &= + \sqrt{-\frac{d^2}{du^2}\log
d_\pm(u)}\cdot\beta(u) + \gamma \alpha(u)
\end{align}
\end{subequations}

The real number $\gamma$ is the  spectral parameter. Let
$\mu_\pm(u) = \sqrt{-\frac{d^2}{du^2}\log d_\pm(u)}$. Each
component
 $\alpha(u)$ and $\beta(u)$ has its own Schrödinger 
equation. These equations are related through a Darboux
transformation, and they take the form:

\begin{subequations}
\begin{align}
- \frac{d^2}{du^2} \alpha(u) + (\mu_\pm(u)^2 -
\mu_\pm'(u))\alpha(u) &= \gamma^2 \alpha(u)\label{eq:schrod1}\\
- \frac{d^2}{du^2} \beta(u) + (\mu_\pm(u)^2 +
\mu_\pm'(u))\beta(u) &= \gamma^2 \beta(u)\label{eq:schrod2}
\end{align}
\end{subequations}

See also \eqref{eq:schrodbis++}, \eqref{eq:schrodbis+-},
\eqref{eq:schrodbis-+}, \eqref{eq:schrodbis--}, for
alternative expressions. The  potentials arising in
\eqref{eq:schrod1}, \eqref{eq:schrod2} are exponentially
increasing as $u\to+\infty$ and exponentially vanishing as
$u\to-\infty$. We have studied the scattering (rather,
reflection) induced by  \eqref{eq:schrod1} and
\eqref{eq:schrod2} on waves originating at $-\infty$ and
being bounced back to $-\infty$ and we have established that
the corresponding phase shifts realize (in the manner
indicated above, and stated more precisely in the French
section) the Fourier cosine and sine transforms as
scattering matrices. A detailed treatment will appear
elsewhere. Here, we shall explain the origin of the systems
\eqref{eq:dirac} and the existence of known solutions for
each value of the spectral parameter: the structure
functions of the Hilbert spaces of entire functions which
had been associated by de~Branges to the Fourier transform
\cite{bra64}. We gave their first explicit formulae in
\cite{jfsonine}  and give here their associated differential
systems \eqref{eq:dirac}: these systems turn out to involve
the determinants $d_\pm(u)$, a fact which thus establishes a
link between the theory of the de~Branges Sonine spaces and
the quantities arising in the study of the random matrices
\cite{mehta}. In \cite{jfhab} we linked through the
co-Poisson formula the study of the Riemann zeta function to
the Sonine spaces. As the reader will perceive, reference
\cite{bra} plays only a distant rôle here. This is despite
the fact that we consider most beautiful the de~Branges
discovery \cite{bra64} of the existence of a natural
isometric expansion for Fourier self-or-skew-reciprocal
functions, and that we are aware of the rich contents of his
general theory \cite{bra}. We are also aware of the rich
relations it holds with other works, such as the
investigations of Krein (see, for example, \cite{dymkean}),
relations which are not that easily perceived from
\cite{bra}. Rather the elucidation of the de~Branges Fourier
Sonine spaces, which was a necessity after \cite{jfcrashp},
was initiated in \cite{jfsonine} and is continued
here, and turns out to provide an opportunity to partially
realize some of the ideas which are internal to our approach
of the problems of the zeta function and which we have
expressed elsewhere \cite{jfprop}, \cite{jfscat},
\cite{jfhab}.

\selectlanguage{french}

\section{Un système différentiel}

Nous continuons la discussion entamée dans la partie en
langue anglaise du système de Dirac \eqref{eq:dirac} et des
équations de Schrödinger \eqref{eq:schrod1},
\eqref{eq:schrod2} associées. Nous utiliserons les notations
suivantes: la transformation de Fourier $\cF$ agit selon
$\cF(f)(y) =
\int_\RR e^{2\pi i xy}f(x)\,dx$. On a $\cF = \cF_+ + i\cF_-$
où $\cF_+$ est la transformation en cosinus et $\cF_-$ est
la transformation en sinus. Soit $a>0$. Nous utiliserons
l'opérateur $P_a$ de restriction à l'intervalle $[-a,a]$,
ainsi que les opérateurs $\cF_a^+ = \cF_+ P_a$, $F_a^+ =
P_a\cF_+ P_a$, $D_a^+ = (F_a^+)^2$, $\cF_a^- = \cF_- P_a$,
$F_a^- =
P_a \cF_- P_a$, $D_a^- = (F_a^-)^2$. Les opérateurs $F_a^+$
et
$D_a^+$ (resp. $F_a^-$ et $D_a^-$) peuvent être vus au choix
comme agissant sur $L^2(-a,a; \frac12 dx)$, sur son
sous-espace pair (resp. impair), ou après un changement
d'échelle sur $L^2(-1,1; \frac12 dx)$ ou encore sur
$L^2(0,1;dx)$. Dans cette dernière représentation le noyau
de l'opérateur $D_a^\pm$ est donné par:
\begin{equation}
\frac{\sin(2\pi a^2 (x-y))}{\pi(x-y)} \pm \frac{\sin(2\pi
a^2 (x+y))}{\pi(x+y)}
\end{equation}

Considérons les uniques
fonctions $\phi_a^+$ et $\phi_a^-$ solutions (disons,
localement intégrables) des équations:

\begin{subequations}\label{eq:phi}
\begin{align}
\phi_a^+ + \cF_a^+ \phi_a^+ &= 2\cos(2\pi ax)\\
\phi_a^- - \cF_a^+ \phi_a^- &= 2\cos(2\pi ax)
\end{align}
\end{subequations}

Ce sont des fonctions paires, entières. Nous mentionnons
également les fonctions impaires qui sont associées à la
transformation en sinus. Elles  vérifient:

\begin{subequations}\label{eq:psi}
\begin{align}
\psi_a^+ + \cF_a^- \psi_a^+ &= 2\sin(2\pi ax)\\
\psi_a^- - \cF_a^- \psi_a^- &= 2\sin(2\pi ax)
\end{align}
\end{subequations}

Pour être spécifique nous considérerons ici la
transformation en cosinus, la discussion étant
essentiellement identique
pour la transformation en sinus. Soit $\delta_x =
x\frac\partial{\partial x} + \frac12$. Un calcul simple
aboutit au système différentiel suivant:

\begin{subequations}\label{eq:phisys}
\begin{align}
a\frac\partial{\partial a} \phi_a^+  = \delta_x \phi_a^- -
(\frac12 + a\phi_a^+(a) + a\phi_a^-(a))\phi_a^+\\
a\frac\partial{\partial a} \phi_a^-  = \delta_x \phi_a^+ -
(\frac12 - a\phi_a^+(a) - a\phi_a^-(a))\phi_a^-
\end{align}
\end{subequations}

La démonstration consiste à montrer que
$a\frac\partial{\partial a} \phi_a^+  -\delta_x \phi_a^-$ et
$a\frac\partial{\partial a} \phi_a^-  - \delta_x \phi_a^+ $
satisfont aussi les équations \eqref{eq:phi} à des facteurs
multiplicatifs près qui sont comme indiqués. La quantité
$a\phi_a^+(a)+ a\phi_a^-(a)$ jouant un rôle important nous
la désignerons par $\mu(a)$:
\begin{equation}
\mu(a) = a\phi_a^+(a)+ a\phi_a^-(a)
\end{equation}
On établit sans peine comme conséquence de
\eqref{eq:phisys}:
\begin{equation}
\frac\partial{\partial a} a \phi_a^+\phi_a^- = \frac12
\frac\partial{\partial x} x ((\phi_a^+)^2 + (\phi_a^-)^2 )
\end{equation}
et à partir de là la formule suivante:

\begin{equation}\label{eq:mu}
\mu(a)^2 = a^2(\phi_a^+(a)+ \phi_a^-(a))^2 = a\frac{d}{da}
a\int_{-a}^a\phi_a^+(x)\phi_a^-(x)\,dx
\end{equation}

\section{Le lien avec les déterminants de Fredholm}

Par une formule connue de la théorie générale des
  déterminants de Fredholm:
\begin{align}
\phi_a^+(a) = +\frac{d}{da} \log\det(1+F_a^+)\ \qquad
\phi_a^-(a) = -\frac{d}{da} \log\det(1-F_a^+)\label{eq:phidet}\\
\psi_a^+(a) = +\frac{d}{da} \log\det(1+F_a^-)\ \qquad
\psi_a^-(a) = -\frac{d}{da} \log\det(1-F_a^-)\label{eq:psidet}
\end{align}
Ainsi, on a aussi:
\begin{equation}
\mu(a) = a\frac{d}{da}
\log\frac{\det(1+F_a^+)}{\det(1-F_a^+)}
\end{equation}

Nous relions maintenant $\mu(a)$ au déterminant de Fredholm
$\det(1-D_a^+)$. Par une formule bien connue on a:
\begin{equation}
\frac{d}{da}\log d_a^+ = -{\rm Tr}((1-
D_a^+)^{-1}\frac{d}{da}D_a^+)
\end{equation}
Or le calcul immédiat du noyau de la dérivée de $D_a^+$
montre qu'il s'agit d'un opérateur de rang un, et donc que
$(1- D_a^+)^{-1}\frac{d}{da}D_a^+$ l'est également. Son 
vecteur propre est donné par la fonction
$(1-D_a^+)^{-1}(2\cos(2\pi a^2x))$ (nous sommes ici sur
$[0,1]$), et on trouve que la valeur propre est égale au
produit scalaire de cette fonction avec $(2a)2\cos(2\pi a^2
x)$. En faisant usage de la factorisation $1-D_a^+ =
(1+F_a^+)(1-F_a^+)$ et du caractère auto-adjoint de ces
opérateurs on obtient finalement la formule:
\begin{equation}
\frac{d}{da} \log d_a^+ = - \int_{-a}^a
\phi_a^+(x)\phi_a^-(x)\,dx
\end{equation}
où l'on est revenu à l'intervalle $[-a,a]$.

La comparaison avec \eqref{eq:mu}, \eqref{eq:phidet},
permet donc d'obtenir le théorème suivant:
\begin{theoreme}
Soit $d_+(u) = d_a^+ = \det(1 - D_a^+)$, et $\mu_+(u) =
\mu(a)$ pour $a = \exp(u)$,
$u\in\RR$. On a:
\begin{align}
\mu_+(u)^2  &= -\frac{d^2}{du^2} \log d_+(u)\\
\mu_+(u)^2 - \mu_+'(u) &= - 2\frac{d^2}{du^2} \log \det(1 +
F_a^+)\label{eq:pot+}\\
\mu_+(u)^2 + \mu_+'(u) &= - 2\frac{d^2}{du^2} \log \det(1 -
F_a^+)\label{eq:pot-}
\end{align}
\end{theoreme}

Dans \eqref{eq:dirac}, on a écrit $\sqrt{-\frac{d^2}{du^2}\log
d_+(u)}$  pour désigner la fonction $u\mapsto \mu_+(u) = \mu(e^u)$.

\section{Le système différentiel des distributions $A_a$ et
$B_a$}

Les fonctions $\phi_a^+$ et $\phi_a^-$ sont aussi des
distributions tempérées dont on peut prendre la transformée
de Fourier. On vérifie de suite que les distributions
tempérées paires:
\begin{subequations}
\begin{align}
A_a &= \frac{\sqrt{a}}2\, (\cF(\phi_a^+) + \phi_a^+ )\\
B_a &= \frac{i\sqrt{a}}2\, (\cF(\phi_a^-) - \phi_a^-)
\end{align}
\end{subequations}
ont la propriété remarquable de s'annuler dans $]-a,a[$ tout
comme leurs transformées de Fourier puisque $A_a$ est
invariante et $B_a$ est anti-invariante. Or un théorème
général \cite[Thm. 1]{jfsonine} nous dit que toute
distribution $D$ (paire ou impaire) ayant cette propriété a
comme transformée de Mellin droite $\int_0^\infty
D(t)t^{-s}dt$ (a priori définie pour $\Re(s)\gg0$) une
fonction entière de $s$ qui a des zéros triviaux en $s =
-2n$, $n\in\NN$ ($-1-2n$, dans le cas impair). Nous noterons
$\cA_a(s)$ et $\cB_a(s)$ les transformées de Mellin de $A_a$
et $B_a$, complétées par le facteur Gamma
$\pi^{-s/2}\Gamma(\frac s2)$. Nous avons démontré dans
\cite{jfsonine} que ces fonctions entières sont en fait les
fonctions de structure (dans une certaine normalisation) de
certains espaces de Hilbert de fonctions entières que
de~Branges a associé à la transformation de Fourier
\cite{bra64}. Avec nos conventions, l'axe de symétrie de ces
espaces est la droite critique, et donc par ce fait même,
les fonction entières $\cA_a$ et $\cB_a$ vérifient
l'hypothèse de Riemann. Nous avons aussi montré que la
densité de leurs zéros est au premier ordre identique à
celle pour les zéros de la fonction dzêta de Riemann
\cite{jftwosys}.

Cela étant dit, à partir du système différentiel
\eqref{eq:phisys} on obtient:
\begin{subequations}
\begin{align}
a\frac\partial{\partial a} A_a  &= i\delta_x B_a -
\mu(a)A_a\\
a\frac\partial{\partial a} B_a  &= -i\delta_x A_a +
\mu(a)B_a
\end{align}
\end{subequations}
Au niveau des transformées de Mellin droites on a $\delta_x =
s-\frac12 = i\gamma$. On obtient donc:
\begin{subequations}
\begin{align}
a\frac\partial{\partial a} \cA_a  &= -\gamma \cB_a -
\mu(a)\cA_a\\
a\frac\partial{\partial a} \cB_a  &= +\gamma \cA_a +
\mu(a)\cB_a
\end{align}
\end{subequations}

Nous poserons $\cA^s(u) = \cA_a(s)$, $\cB^s(u) = \cB_a(s)$,
$\mu_+(u) = \mu(a)$ avec $a=\exp(u)$, $s = \frac12 +
i\gamma$. Ainsi:

\begin{theoreme}
Les fonctions $\cA^s(u)$ et $\cB^s(u)$ de la variable
$u\in\RR$ vérifient le système différentiel suivant
($s=\frac12+i\gamma$):
\begin{subequations}
\begin{align}
\frac{d}{du} \cA^s(u)  &= -\gamma \cB^s(u) -
\mu_+(u)\cA^s(u) \\
\frac{d}{du} \cB^s(u) &= +\gamma \cA^s(u) + \mu_+(u)\cB^s(u)
\end{align}
\end{subequations}
\end{theoreme}

Les équations de Schrödinger associées en découlent:
\begin{subequations}
\begin{align}
- \frac{d^2}{du^2}  \cA^s(u) + \Big(- 2\frac{d^2}{du^2} \log \det(1 +
F_a^+)\Big) \cA^s(u) &= \gamma^2  \cA^s(u)\label{eq:schrodbis++}\\
- \frac{d^2}{du^2}  \cB^s(u) + \Big(- 2\frac{d^2}{du^2} \log \det(1 -
F_a^+)\Big) \cB^s(u)&= \gamma^2 \cB^s(u)\label{eq:schrodbis+-}
\end{align}
\end{subequations}

Pour la transformation en sinus on obtiendrait:
\begin{subequations}
\begin{align}
- \frac{d^2}{du^2}  \alpha(u) + \Big(- 2\frac{d^2}{du^2} \log \det(1 +
F_a^-)\Big) \alpha(u) &= \gamma^2  \alpha(u)\label{eq:schrodbis-+}\\
- \frac{d^2}{du^2}  \beta(u) + \Big(- 2\frac{d^2}{du^2} \log \det(1 -
F_a^-)\Big) \beta(u)&= \gamma^2 \beta(u)\label{eq:schrodbis--}
\end{align}
\end{subequations}

On
notera que  si l'on normalise les fonctions de structure
$\cA$ et $\cB$ des espaces de de~Branges en imposant
$\cA_a(\frac12) = 1$ alors les équations du deuxième ordre
vérifiées par $\cA$ et $\cB$ prennent la forme d'équations
de Sturm-Liouville.

\section{Compléments}

La théorie générale de de~Branges \cite{bra} établit
l'existence d'équations intégrales (qui peuvent parfois
prendre, comme il était connu a priori que cela serait le
cas ici \cite{bra}, une forme différentielle) pour certaines
chaînes de sous-espaces (dont les éléments sont aussi des
fonctions entières) d'un même espace de Hilbert
$L^2(\RR,d\Delta(\gamma))$. Elle les utilise pour obtenir
une représentation isométrique de
$L^2(\RR,d\Delta(\gamma))$, apportant ainsi une vaste
généralisation au théorème de Paley-Wiener. Dans le cas
particulier qui est considéré ici nous avons vérifié que
l'expansion isométrique des espaces de Sonine (associés à la
transformation en cosinus) prend la forme suivante:
l'application
\begin{equation}\label{eq:iso}
\begin{bmatrix}
\alpha(u)\\ \beta(u)
\end{bmatrix}
\mapsto
f(s) = \int_{-\infty}^{+\infty}
\Big(\cA^s(u)\alpha(u)+\cB^s(u)\beta(u)\Big)\,du
\end{equation}
est isométrique de $L^2(\RR\to\RR^2;\frac{du}2)$ sur
$L^2(\Re(s)=\frac12;d\Delta(\gamma))$ avec: 
\begin{equation}
d\Delta(\gamma) = \frac1{\left|\pi^{-\frac s2}\Gamma(\frac
s2)\right|^2}\frac{|ds|}{2\pi}, \qquad  s = \frac12 +
i\gamma.
\end{equation}
Ce type d'expansion isométrique est aussi intimement lié à
des considérations antérieures de Krein, comme il est
expliqué dans le livre de Dym et McKean
\cite{dymkean}. 

Les potentiels pour les équations de Schrödinger sont
exponentiellement croissants en $u\to+\infty$ et s'annulent
exponentiellement en $u\to-\infty$: cela suggère, au-delà
des considérations liées à l'expansion isométrique,
d'aborder le problème du scattering d'ondes provenant de
$-\infty$ et y étant renvoyées. Krein a  utilisé ses
techniques pour aborder sur un demi-axe les problèmes de scattering
également traités par les approches de Gel'fand-Levitan et
Marchenko, dans un cadre où la mesure spectrale
$d\Delta(\gamma)$ satisfait la condition
$\int\frac{d\Delta(\gamma)}{1+\gamma^2}<\infty$, ce qui
n'est pas le cas ici, bien qu'un grand nombre de
considérations restent valables. En fait la différence
essentielle avec la situation habituellement envisagée est
qu'ici, heuristiquement, plus l'onde est énergétique, plus
profondément elle pénètrera dans le potentiel et plus
longtemps elle mettra pour en être renvoyée: de sorte que
l'on s'attend à ce que les décalages de phases soient
divergents lorsque la fréquence $\gamma$ tend vers
$\pm\infty$.

Ce problème, pour \eqref{eq:schrodbis++} par exemple,  est
rendu précis par le fait que la condition ``être de carré
intégrable en $+\infty$'' est une condition aux bords
appropriée en $+\infty$, avec pour unique solution (à une
constante multiplicative près) de fréquence réelle $\gamma$
la fonction $u\mapsto \cA^s(u)$. Le calcul de la matrice de
scattering découle alors de la détermination, pour chacune
des équations de Schrödinger \eqref{eq:schrodbis++},
\eqref{eq:schrodbis+-}, \eqref{eq:schrodbis-+},
\eqref{eq:schrodbis--}, des solutions satisfaisant la
``condition de Jost'' en $-\infty$. Par exemple, on trouve,
pour l'équation  \eqref{eq:schrodbis++}, que la solution
(complexe) qui se comporte asymptotiquement en $-\infty$
comme $\exp(-i\gamma u)$ est donnée exactement par
($a=\exp(u)$, $s=\frac12 + i\gamma$, $\Re(s)<1$):
\begin{equation}\label{eq:jost}
u\mapsto e^{-i\gamma u} - e^{u/2} \int_0^a \phi_a^+(x)x^{-s}\,dx
\end{equation}
et à partir de là on établit que la ``matrice de
scattering'' de $-\infty$ vers $-\infty$ est, pour
l'équation \eqref{eq:schrodbis++}, exactement donnée par la
fonction $ \pi^{-(1-s)/2}\Gamma((1-s)/2) /
\pi^{-s/2}\Gamma(s/2)=\chi_+(s)$, $s=\frac12 + i\gamma$, qui
donne la représentation multiplicative (sur la droite
critique) de la transformation en cosinus.  Pour l'équation
\eqref{eq:schrodbis+-} on trouve $-\chi_+(\frac12+i\gamma)$
comme matrice de scattering. Pour les équations
\eqref{eq:schrodbis-+} et \eqref{eq:schrodbis--} associées à
la transformation en sinus on remplace $\chi_+$ par la
fonction $\chi_-$ donnant la représentation multiplicative
de la transformation en sinus. Peut-être existe-t-il un moyen
d'utiliser cette information pour remonter aux potentiels
eux-mêmes et obtenir ainsi peut-être une autre approche
(spécifique au cas d'un unique intervalle $[-a,a]$) aux
résultats de Deift-Its-Zhou \cite{diz} et des autres auteurs
cités.

Ces ``matrices de scattering'' sont aussi les fonctions qui
apparaissent  dans les équations fonctionnelles de la
fonction dzêta de Riemann et des séries $L$ de
Dirichlet. Signalons à ce sujet que la présente Note est la
continuation logique de nos travaux antérieurs
\cite{jfprop}, \cite{jfcond} (voir \cite{jfcrasexp} pour un
résumé) sur l'opérateur conducteur $\log|x|+\log|y|$,
opérateur qui entre autres choses, révéla que la
transformation de Fourier était une symétrie des formules
explicites.

\bibliographystyle{amsplain}

\end{document}